\documentclass[final,a4paper,11pt,reqno,twoside]{amsart}
\usepackage{riro}
\begin{document} 
\begin{merci}%
 The first author is financed by the ANR project ``Aspects Arithm\'etiques des Matrices Al\'eatoires et du Chaos Quantique''. He would like to thank the University of Nottingham, where this work has been finished, for its hospitality. The second one is supported by the ANR project ``Modunombres''. Both authors would like to thank the anonymous referee of \cite{RiRo} for suggesting them this problem. %
\end{merci}
\section{Introduction and statement of the result}%
\subsection{Description of the families of $L$-functions studied}%
The purpose of this paper is to compute the lower order terms of some particular statistics associated to low-lying zeros of several families of symmetric power $L$-functions in the level aspect: the one-level densities. First of all, we give a short description of these families. To any primitive holomorphic cusp form $f$ of prime level $q$ and even weight\footnote{In this paper, the weight $\kappa$ is a \emph{fixed} even integer and the level $q$ goes to infinity among the prime numbers.} $\kappa\geq 2$ (see \cite[\S~2.1]{RiRo} for the automorphic background) say $f\in\prim{\kappa}{q}$, one can associate its $r$-th symmetric power $L$-function denoted by $L(\sym^rf,s)$ for any integer $r\geq 1$. It is given by the following absolutely convergent and non-vanishing Euler product of degree $r+1$ on $\Re{s}>1$ %
\[%
L(\sym^rf,s)=\prod_{p\in\prem}L_p(\sym^rf,s) %
\]
where %
\[%
L_p(\sym^rf,s)=\prod_{i=0}^r\left(1-\frac{\alpha_f(p)^{i}\beta_f(p)^{r-i}}{p^{s}}\right)^{-1} %
\]
for any prime number $p$. From now on, $\alpha_f(p)$, $\beta_f(p)$ are the Satake parameters of $f$ at the prime number $p$ and $\left(\lambda_f(n)\right)_{n\geq 1}$ is its sequence of Hecke eigenvalues, which is arithmetically normalised: $\lambda_f(1)=1$ and $\abs{\lambda_f(p)}\leq 2$ for any prime $p$. We also define \cite[(3.16) and (3.17)]{CoMi} a local factor at $\infty$ which is given by a product of $r+1$ Gamma factors namely %
\begin{equation}\label{eq_r_odd}%
L_\infty(\sym^rf,s)=\prod_{0\leq a\leq(r-1)/2}\fGamma_{\R}\left(s+(2a+1)(\kappa-1)/2\right)\fGamma_{\R}\left(s+1+(2a+1)(\kappa-1)/2\right) %
\end{equation}
if $r$ is odd and %
\begin{equation}\label{eq_r_even}%
L_\infty(\sym^rf,s)=\fGamma_{\R}(s+\mu_{\kappa,r})\prod_{1\leq a\leq r/2}\fGamma_{\R}\left(s+a(\kappa-1)\right)\fGamma_{\R}\left(s+1+a(\kappa-1)\right) %
\end{equation}
if $r$ is even where %
\[%
\mu_{\kappa,r}=\begin{cases}%
1 & \text{if } r(\kappa-1)/2 \text{ is odd,} \\%
0 & \text{otherwise.}\\%
\end{cases}
\]
The completed $L$-function is defined by %
\begin{equation*}%
\Lambda(\sym^rf,s)=\left(q^{r}\right)^{s/2}L_\infty(\sym^rf,s)L(\sym^rf,s)%
\end{equation*}
and $q^r$ is the arithmetic conductor. We will need some control on the analytic behaviour of this function. Unfortunately, such information is not currently known in all generality. We sum up our main assumption in the following statement. %
\begin{hypothesis}\label{hypohypo} %
The function $\Lambda\left(\sym^rf,s\right)$ is a \emph{completed $L$-function} in the sense that it satisfies the following \emph{nice} analytic properties: %
\begin{itemize}%
\item it can be extended to an holomorphic function of order $1$ on $\C$, %
\item it satisfies a functional equation of the shape %
\[%
\Lambda(\sym^rf,s)=\epsilon\left(\sym^rf\right)\Lambda(\sym^rf,1-s) %
\]
where the sign $\epsilon\left(\sym^rf\right)=\pm 1$ of the functional equation is given by %
\begin{equation}\label{valueofsign}%
\epsilon\left(\sym^rf\right)=%
\begin{cases}
+1 & \text{if $r$ is even},\\
\epsilon_f(q)\times\epsilon(\kappa,r) & \text{otherwise}%
\end{cases}
\end{equation}
with %
\[%
\epsilon(\kappa,r)= i^{\left(\frac{r+1}{2}\right)^2(\kappa-1)+\frac{r+1}{2}}=%
\begin{cases}
i^{\kappa} & \text{if $\;r\equiv 1\pmod{8}$,} \\ %
-1 & \text{if $\;r\equiv 3\pmod{8}$}, \\%
-i^{\kappa} & \text{if $\;r\equiv 5\pmod{8}$}, \\%
+1 & \text{if $r\;\equiv 7\pmod{8}$}%
\end{cases}
\]
and $\epsilon_f(q)=-\sqrt{q}\lambda_f(q)=\pm 1$. %
\end{itemize}
\end{hypothesis}
\begin{remint}%
Hypothesis $\Nice(r,f)$ is known for $r=1$ (E.~Hecke \cite{MR1513069,MR1513122,MR1513142}), $r=2$ thanks to the work of S.~Gelbart and H.~Jacquet \cite{GeJa} and $r=3,4$ from the works of H.~Kim and F.~Shahidi \cite{KiSh1,KiSh2,Ki}. %
\end{remint}
We aim at studying the lower order terms of the one-level density for the family of $L$-functions given by %
\begin{equation*}
\bigcup_{\text{$q$ prime}}\left\{L(\sym^rf,s), f\in\prim{\kappa}{q}\right\} %
\end{equation*}
for any integer $r\geq 1$. %

\subsection{One-level densities of these families}%
The purpose of this work is to determine the lower order terms of the one-level densities associated to these families of $L$-functions. Let us give the statement of our result, in which $\nu$ is a positive real number, $\Phi$ is an even Schwartz function, whose Fourier transform $\widehat{\Phi}$ is compactly supported in $[-\nu,+\nu]$ (denoted by $\Phi\in\Schwartz_\nu(\R)$) and $f$ is a primitive holomorphic cusp form of prime level $q$ and even weight $\kappa\geq 2$ for which hypothesis $\Nice(r,f)$ holds\footnote{Note that we do not assume any Generalised Riemann Hypothesis for the symmetric power $L$-functions. }. We refer to \cite[\S~2.2]{RiRo} for the probabilistic background. Note that, thanks to Fourier inversion formula, such a function $\Phi$ can be extended to an entire even function  which satisfies %
\begin{equation}\label{estim}%
\forall s\in\C,\quad\Phi(s)\ll_n\frac{\exp{(\nu\abs{\Im{s}})}}{(1+\abs{s})^n} %
\end{equation}
for any integer $n\geq 0$. The \emph{one-level density} (relatively to $\Phi$) of $\sym^rf$ is defined by %
\begin{equation*}%
D_{1,q}[\Phi;r](f)=\sum_{\rho,\;\Lambda(\sym^rf,\rho)=0}\Phi\left(\frac{\log{\left(q^r\right)}}{2i\pi}\left(\Re{\rho}-\frac{1}{2}+i\Im{\rho}\right)\right)%
\end{equation*}
where the sum is over the non-trivial zeros $\rho$ of $L(\sym^rf,s)$ repeated with multiplicities. The \emph{asymptotic expectation} of the one-level density is by definition %
\begin{equation*}%
\lim_{\substack{q\;\text{prime} \\ q\to+\infty}}{\sum_{f\in\prim{\kappa}{q}}}\omega_q(f)D_{1,q}[\Phi;r](f)
\end{equation*}
where $\omega_q(f)=\frac{\fGamma(\kappa-1)}{(4\pi)^{\kappa-1}\scal{f}{f}_q}$ is the harmonic weight of $f$. Before stating our result, let us define the following constants: %
\begin{align}%
C_{\text{PNT}} & =\left(1+\int_1^{+\infty}\frac{\theta(t)-t}{t^2}\dd t\right), \label{eq_CPNT} \\
C & =\sum_{p\in\mathcal{P}}\frac{\log{p}}{p^{3/2}-p}, \label{eq_Cr} \\
C_{\infty} &  =-(r+1)\log{\pi}+C_{\fGamma} \label{eq_Cinfty}
\end{align}
where $\theta$ is the first Chebyshev function: %
\[%
 \theta(t)= \sum_{\substack{\text{$p$ prime}\\ p\leq t}}\log p, %
\]
\begin{equation}\label{eq_C_Gamma_odd}%
C_{\fGamma}=\sum_{0\leq a\leq(r-1)/2}\left\{\left(\frac{\fGamma'}{\fGamma}\right)\left(\frac{1}{4}+\frac{(2a+1)(\kappa-1)}{4}\right)+\left(\frac{\fGamma'}{\fGamma}\right)\left(\frac{1}{4}+\frac{1}{2}+\frac{(2a+1)(\kappa-1)}{4}\right)\right\}
\end{equation}
if $r$ is odd and %
\begin{equation}\label{eq_C_Gamma_even}%
C_{\fGamma}=\left(\frac{\fGamma'}{\fGamma}\right)\left(\frac{1}{4}+\frac{\mu_{\kappa,r}}{2}\right)+\sum_{1\leq a\leq r/2}\left\{\left(\frac{\fGamma'}{\fGamma}\right)\left(\frac{1}{4}+\frac{a(\kappa-1)}{2}\right)+\left(\frac{\fGamma'}{\fGamma}\right)\left(\frac{1}{4}+\frac{1}{2}+\frac{a(\kappa-1)}{2}\right)\right\} %
\end{equation}
if $r$ is even. %
\begin{theoint}\label{thm_B}%
Let $r\geq 1$ be any integer and $\epsilon=\pm 1$. We assume that hypothesis $\Nice(r,f)$ holds for any prime number $q$ and any primitive holomorphic cusp form of level $q$ and even weight $\kappa\geq 2$. Let%
\[%
\nu_{1,\mathrm{max}}(r,\kappa,\theta_0)=\left(1-\frac{1}{2(\kappa-2\theta_0)}\right)\frac{2}{r^2}%
\]
with $\theta_0=7/64$. If $\nu<\nu_{1,\mathrm{max}}(r,\kappa,\theta_0)$ then the asymptotic expectation of the one-level density is %
\[%
\left[\widehat{\Phi}(0)+\frac{(-1)^{r+1}}{2}\Phi(0)\right]+\left[C_{\infty}-2(-1)^rC_{\text{PNT}}-2\delta_{2\mid r}C\right]\frac{\widehat{\Phi}(0)}{\log q^r}+O\left(\frac{1}{\log^3(q^r)}\right). %
\]
\end{theoint}
\begin{remint}
The main terms of the asymptotic expectation of these one-level densities have already been found in \cite{RiRo} (see Theorem B). The new information is the lower order terms namely terms of size $1/\log{(q^r)}$. %
\end{remint}
\begin{remint}
Note that $\theta_0=7/64$ is the best known approximation towards Ramanujan-Peterson-Selberg's conjecture (see \cite[hypothesis $\Hy_2(\theta)$ page 16]{RiRo}) thanks to the works of H. Kim, F. Shahidi and P. Sarnak (\cite{KiSh2,Ki}). The value $\theta=0$ is expected. %
\end{remint}
\begin{remint}
It is clear from the proof of Theorem \ref{thm_B} that the same result holds for the signed families with the same restriction on the support as in \cite{RiRo}. %
\end{remint}
\begin{remint}
The particular case $r=1$ has already been investigated by S.J.~Miller\cite{Mil}. %
\end{remint}

\begin{notations}%
We write $\prem$ for the set of prime numbers and the main parameter in this paper is a prime number $q$, whose name is the level, which goes to infinity among $\prem$. Thus, if $f$ and $g$ are some $\C$-valued functions of the real variable then the notations $f(q)\ll_{A}g(q)$ or $f(q)=O_A(g(q))$ mean that $\abs{f(q)}$ is smaller than a ``constant'' which only depends on $A$ times $g(q)$ at least for $q$ a large enough prime number.
\end{notations}

\section{Chebyshev polynomials and Hecke eigenvalues}\label{cheby} 

Recall that the general facts about holomorphic cusp forms can be found in \cite[\S~2.1]{RiRo}. Let $p\neq q$ a prime number and $f\in\prim{\kappa}{q}$. Denote by $\chi_{\St}$ the character of the standard representation $\St$ of $\sud$. By the work of Deligne, there exists $\theta_{f,p}\in[0,\pi]$ such that %
\[%
 \lambda_f(p)=\chi_{\St}\begin{pmatrix}e^{i\theta_{f,p}} & 0\\ 0 & e^{-i\theta_{f,p}}\end{pmatrix}. %
\]
Moreover the multiplicativity relation reads %
\begin{equation}\label{eq_ams}%
 \lambda_{f}(p^\nu)=\chi_{\sym^\nu}\begin{pmatrix}e^{i\theta_{f,p}} & 0\\ 0 & e^{-i\theta_{f,p}}\end{pmatrix}=X_\nu\left(\chi_{\St}\begin{pmatrix}e^{i\theta_{f,p}} & 0\\ 0 & e^{-i\theta_{f,p}}\end{pmatrix}\right)=X_\nu\left(\lambda_f(p)\right) %
\end{equation}
where $\chi_{\sym^\nu}$ is the character of the irreducible representation $\sym^\nu\St$ of $\sud$ and the polynomials $X_\nu$ are defined by their generating series %
\begin{equation}\label{eq_genche}%
\sum_{\nu\geq 0}X_\nu(x)t^\nu=\frac{1}{1-xt+t^2}.
\end{equation}
They are equivalentely defined by %
\begin{equation}\label{eq_che}%
X_\nu(2\cos\theta)=\frac{\sin{((\nu+1)\theta)}}{\sin{(\theta)}}. %
\end{equation}  
These polynomials are known as \emph{Chebyshev polynomials of second kind}. Each $X_\nu$ has degree $\nu$, is even if $\nu$ is even and odd otherwise. The family $\left(X_\nu\right)_{\nu\geq 0}$ is a basis for the polynomial vector space $\Q[T]$, orthonormal with respect to the inner product %
\[%
\scalst{P}{Q}=\frac{1}{\pi}\int_{-2}^{2}P(x)Q(x)\sqrt{1-\frac{x^2}{4}}\dd x. %
\] 
The following proposition lists Chebyshev polynomials' needed properties for this work. %
\begin{proposition}\label{prop_cheby}\strut\newline%
\begin{itemize}
\item
If $\varpi\geq 0$ is any integer then %
\begin{equation}\label{eq_lintch}%
X_r^\varpi=\sum_{j=0}^{r\varpi}x(\varpi,r,j)X_j %
\end{equation}
with %
\begin{equation}\label{eq_valx}%
x(\varpi,r,j)= \scalst{X_r^\varpi}{X_j}=\frac{2}{\pi}\int_0^\pi\frac{\sin^\varpi{((r+1)\theta)}\sin{((j+1)\theta)}}{\sin^{\varpi-1}{(\theta)}}\dd\theta.%
\end{equation}
In particular, %
\begin{equation}\label{propriox}%
x(\varpi,r,j)=\begin{cases}%
0 & \text{if $j\equiv r\varpi+1\pmod{2}$}, \\%
\frac{\binom{\varpi}{\varpi/2}}{1+\varpi/2} & \text{if $\varpi$ is even, $r=1$ and $j=0$.} \\%
\end{cases}%
\end{equation}
\item %
If $\alpha$ is a complex number of norm $1$ and $n\geq 0$ is an integer then %
\begin{equation}\label{eq_S4}%
\alpha^n+\alpha^{-n}=\begin{cases}%
2X_0(\alpha+\alpha^{-1}) & \text{if $n=0$,} \\%
X_1(\alpha+\alpha^{-1}) & \text{if $n=1$,} \\%
X_n(\alpha+\alpha^{-1})-X_{n-2}(\alpha+\alpha^{-1}) & \text{otherwise.}%
\end{cases}
\end{equation}
\item %
If $\alpha$ is a complex number of norm $1$ and $r, n\geq 1$ are some integers then %
\begin{align}\label{eq_S}%
S(\alpha;n,r)=\sum_{j=0}^r\alpha^{n(2j-r)} 
&= \delta_{2\mid r}+\sum_{\substack{1\leq j\leq r  \\%
j\equiv r\pmod{2}}}\left[\alpha^{jn}+\alpha^{-jn}\right]  \\%
&= \sum_{\substack{0\leq j\leq r \\%
j\equiv r\pmod{2}}}\left[X_{jn}(\alpha+\alpha^{-1})-X_{jn-2}(\alpha+\alpha^{-1})\right] \\%
&= X_r(\alpha^n+\alpha^{-n})  %
\end{align}
where $X_{-1}=X_{-2}=0$ by convention. %
\item%
If $r\geq 1$ and $n\geq 1$ are some integers then %
\begin{equation}\label{eq_S3}%
\sum_{\substack{0\leq j\leq r \\%
j\equiv r\pmod{2}}}\left[X_{jn}-X_{jn-2}\right]=\sum_{j=0}^r(-1)^jX_{n-2}^jX_{n(r-j)}%
\end{equation}%
where $X_{-1}=X_{-2}=0$ by convention. %
\item%
If $\ell\geq 0$ is an integer then %
\begin{equation}\label{eq_lin}%
X_\ell=\sum_{\substack{0\leq u\leq\ell \\%
u\equiv\ell\pmod{2}}}(-1)^{(\ell-u)/2}\binom{(\ell+u)/2}{u}T^u. %
\end{equation}
\end{itemize}
\end{proposition}
\begin{proof}[\proofname{} of proposition \ref{prop_cheby}]%
The first point follows from the fact that $X_r^\varpi$ is an polynomial of degree $r\varpi$, which is even if $r\varpi$ is even and odd otherwise. Thus, \eqref{eq_lintch} is the expansion of this polynomial in the orthonormal basis $\left(X_j\right)_{0\leq j\leq r\varpi}$. The second point follows from the equality %
\begin{equation*}
2\cos{(n\theta)}\sin{(\theta)}=\sin{((n+1)\theta)}-\sin{((n-1)\theta)}. %
\end{equation*}
If $\alpha=\exp{(i\theta)}$ then this equality combined with \eqref{eq_che} lead to %
\begin{equation*}
2\cos{(n\theta)}=X_{n}(2\cos{\theta})-X_{n-2}(2\cos{\theta}), %
\end{equation*}
which is the desired result since $2\cos{\theta}=\alpha+\alpha^{-1}$ and $2\cos{(n\theta)}=\alpha^{n}+\alpha^{-n}$. The third point is a direct consequence of the second one, of the direct computation %
\begin{equation*}%
S(\alpha;n,r)=\frac{\alpha^{n(r+1)}-\alpha^{-n(r+1)}}{\alpha^{n}-\alpha^{-n}} %
\end{equation*}
and of %
\[%
X_r(\alpha^n+\alpha^{-n})  =  X_r(2\cos{(n\theta)}) =  \frac{\alpha^{n(r+1)}-\alpha^{-n(r+1)}}{\alpha^{n}-\alpha^{-n}}%
\]%
if $\alpha=\exp{(i\theta)}$. The fourth point is easily deduced from the fact that %
\begin{equation*}%
S(\alpha;n,r)=\sum_{j=0}^r(-1)^jX_{n-2}^j(\alpha+\alpha^{-1})X_{n(r-j)}(\alpha+\alpha^{-1}) %
\end{equation*}
for any complex number $\alpha$ of norm $1$. Let us prove the previous equality. According to \cite[Page 727, first and second equations]{MR2388795}, %
\begin{equation*}%
\sum_{r\geq 0}X_{nr}(\alpha+\alpha^{-1})t^r=\left[1+X_{n-2}(\alpha+\alpha^{-1})t\right]\sum_{r\geq 0}X_r(\alpha^n+\alpha^{-n})t^r. %
\end{equation*}
As a consequence, %
\begin{equation*}
X_{nr}(\alpha+\alpha^{-1})=X_r(\alpha^n+\alpha^{-n})+X_{n-2}(\alpha+\alpha^{-1})X_{r-1}(\alpha^n+\alpha^{-n}), %
\end{equation*}
which implies %
\begin{equation*}
X_r(\alpha^n+\alpha^{-n})=\sum_{j=0}^r(-1)^jX_{n-2}^j(\alpha+\alpha^{-1})X_{n(r-j)}(\alpha+\alpha^{-1}). %
\end{equation*}
The last point is obtained by developping \eqref{eq_genche} as an entire series in $x$. %
\end{proof}

\section{Riemann's explicit formula for symmetric power $L$-functions}\label{sec_explicit}%

To study $D_{1,q}[\Phi;r](f)$ for any $\Phi\in\Schwartz_\nu(\R)$, we transform this sum over zeros into a sum over primes in the next proposition. In other words, we establish an explicit formula for symmetric power $L$-functions. %
\begin{proposition}\label{explicit}%
Let $r\geq 1$ and $f\in \prim{\kappa}{q}$ for which hypothesis $\Nice(r,f)$ holds and let $\Phi\in\Schwartz_\nu(\R)$. We have %
\begin{multline*}
D_{1,q}[\Phi;r](f)=%
\left[\widehat{\Phi}(0)+\frac{(-1)^{r+1}}{2}\Phi(0)\right]+\frac{\widehat{\Phi}(0)}{\log{(q^r)}}\left[C_{\infty}+2(-1)^{r+1}C_{\text{PNT}}-2\delta_{2\mid r}C\right] \\
+P_q^1[\Phi;r](f)+\sum_{m=0}^{r-1}(-1)^mP_q^2[\Phi;r,m](f)+P_q^3[\Phi;r](f)+O\left(\frac{1}{\log^3{\left(q^r\right)}}\right)
\end{multline*}
where $C_{\text{PNT}}$ is defined in \eqref{eq_CPNT}, $C$ in \eqref{eq_Cr}, $C_\infty$ in \eqref{eq_Cinfty} whereas %
\begin{align*}%
P_{q}^1[\Phi;r](f) & = %
-\frac{2}{\log{\left(q^r\right)}}\sum_{\substack{p\in\prem \\ p\nmid q}}%
\lambda_{f}\left(p^r\right)\frac{\log{p}}{\sqrt{p}}\widehat{\Phi}\left(\frac{\log{p}}{\log{\left(q^r\right)}}\right), \\
P_q^2[\Phi;r,m](f) & =  %
-\frac{2}{\log{\left(q^r\right)}}\sum_{\substack{p\in\prem \\ p\nmid q}}%
\lambda_f\left(p^{2(r-m)}\right)\frac{\log{p}}{p}\widehat{\Phi}\left(\frac{2\log{p}}{\log{\left(q^r\right)}}\right) \\
P_q^3[\Phi;r](f) & =  %
-\frac{2}{\log{\left(q^r\right)}}\sum_{\substack{p\in\prem \\ p\nmid q}}%
\sum_{n\geq 3}\left[\sum_{\substack{1\leq j\leq r \\
j\equiv r\pmod{2}}}\left(\lambda_f(p^{jn})-\lambda_f(p^{jn-2})\right)\right]\frac{\log{p}}{p^{n/2}}\widehat{\Phi}\left(\frac{n\log{p}}{\log{\left(q^r\right)}}\right)
\end{align*}
for any integer $m\in\{0,\ldots,r-1\}$.
\end{proposition}
\begin{proof}[\proofname{} of proposition \ref{explicit}] Let %
\[%
G(s)=\Phi\left(\frac{\log{\left(q^r\right)}}{2i\pi}\left(s-\frac{1}{2}\right)\right). %
\]
From \cite[eq. (4.11) and (4.14)]{IwLuSa} we get %
\begin{multline}\label{eq_etape}
 D_{1,q}[\Phi;r](f)=\widehat{\Phi}(0)-(r+1)\frac{\log\pi}{\log q^r}\widehat{\Phi}(0)\\-\frac{2}{\log q^r}\sum_{p\in\prem}\sum_{m=1}^{+\infty}\left[\sum_{j=0}^r\alpha_f(p)^{jm}\beta_f(p)^{(r-j)m}\right]\widehat{\Phi}\left(\frac{m\log p}{\log q^r}\right)\frac{\log p}{p^{m/2}}\\+\frac{\widehat{\Phi}(0)}{\log q^r}\sum_{j=0}^r\frac{\Gamma'}{\Gamma}\left(\frac{1}{4}+\frac{\mu_j}{2}\right)+O\left(\frac{1}{\log^3q}\right). %
\end{multline}
Let us focus on the \underline{third term in\eqref{eq_etape}}. Not that the contribution of the prime $q$ is given by %
\[%
 -\frac{2}{r}\sum_{m=1}^{+\infty}\left(\frac{\lambda_f(q)^r}{\sqrt{q}}\right)^m\widehat{\Phi}\left(\frac{m}{r}\right)\ll\frac{1}{q^{(r+1)/2}} %
\]
and for $p\neq q$ we use %
\[%
 \sum_{j=0}^r\alpha_f(p)^{jm}\beta_f(p)^{(r-j)m}=S\left(\alpha_f(p);m,r\right) %
\]
with the notation of \eqref{eq_S}. We obtain %
\[%
S\left(\alpha_f(p);1,r\right)= X_r\left(\alpha_f(p)+\alpha_f(p)^{-1}\right)=\lambda_f\left(p^r\right)
\]%
according to \eqref{eq_ams} and 
\begin{align*}
S\left(\alpha_f(p);2,r\right) &= \sum_{\substack{0\leq j\leq r \\ j\equiv r\pmod{2}}}X_{2j}\left(\alpha_f(p)+\alpha_f(p)^{-1}\right)-X_{2j-2}\left(\alpha_f(p)+\alpha_f(p)^{-1}\right) \\
&=  \sum_{j=0}^r(-1)^jX_{2(r-j)}\left(\alpha_f(p)+\alpha_f(p)^{-1}\right) \quad(\text{\emph{cf.}~\eqref{eq_S3}}) \\
&=  \sum_{m=0}^{r-1}(-1)^m\lambda_f\left(p^{2(r-m)}\right)+(-1)^r. %
\end{align*}
As a consequence, %
\begin{multline}\label{etape2}
\sum_{\substack{p\in\prem\\ p\neq q}}\sum_{m=1}^{+\infty}\left[\sum_{j=0}^r\alpha_f(p)^{jm}\beta_f(p)^{(r-j)m}\right]=\sum_{\substack{p\in\mathcal{P}\\ p\nmid q}}\frac{\lambda_f(p^r)\log{p}}{p^{1/2}}\widehat{\Phi}\left(\frac{\log{p}}{\log{\left(q^r\right)}}\right) \\ +\sum_{\substack{p\in\mathcal{P}\\ p\nmid q}}\left(\sum_{m=0}^{r-1}(-1)^m\lambda_f\left(p^{2(r-m)}\right)\right)\frac{\log{p}}{p}\widehat{\Phi}\left(\frac{\log{(p^2)}}{\log{\left(q^r\right)}}\right) \\ +(-1)^r\sum_{\substack{p\in\prem \\ p\nmid q}}\frac{\log{p}}{p}\widehat{\Phi}\left(\frac{\log{(p^2)}}{\log{\left(q^r\right)}}\right) \\ +\sum_{\substack{p\in\mathcal{P}\\ p\nmid q}}\sum_{n\geq 3}S\left(\alpha_f(p);n,r\right)\frac{\log{p}}{p^{n/2}}\widehat{\Phi}\left(\frac{\log{(p^n)}}{\log{\left(q^r\right)}}\right).
\end{multline}
We have isolated the three first terms in \eqref{etape2} since they may contribute as main terms and not only lower order terms. Let us estimate the \underline{third term of \eqref{etape2}}. By partial summation, this term equals, up to $O(q^{-0.9})$,
 \begin{equation*}
(-1)^r\int_1^{+\infty}\frac{\theta(t)}{t^2}\left(\widehat{\Phi}\left(\frac{2\log{t}}{\log{\left(q^r\right)}}\right)-\frac{2}{\log{\left(q^r\right)}}\widehat{\Phi}^\prime\left(\frac{2\log{t}}{\log{\left(q^r\right)}}\right)\right)\dd t:=S_3.
\end{equation*}
Then, %
\begin{multline*}
S_3 =(-1)^r\int_1^{+\infty}\left(\widehat{\Phi}\left(\frac{2\log{t}}{\log{\left(q^r\right)}}\right)-\frac{2}{\log{\left(q^r\right)}}\widehat{\Phi}^\prime\left(\frac{2\log{t}}{\log{\left(q^r\right)}}\right)\right)\frac{\dd t}{t} \\ +(-1)^r\int_1^{+\infty}\frac{\theta(t)-t}{t}\left(\widehat{\Phi}\left(\frac{2\log{t}}{\log{\left(q^r\right)}}\right)-\frac{2}{\log{\left(q^r\right)}}\widehat{\Phi}^\prime\left(\frac{2\log{t}}{\log{\left(q^r\right)}}\right)\right)\frac{\dd t}{t}.
\end{multline*}
Since $\widehat{\Phi}(u)=\widehat{\Phi}(0)+O(u^2)$ and $\widehat{\Phi}^\prime(u)\ll\vert u\vert$, we get
\begin{multline*}
S_3=(-1)^r\frac{\log{\left(q^r\right)}}{2}\int_0^{+\infty}\widehat{\Phi}(u)\dd u-(-1)^r\int_0^{+\infty}\widehat{\Phi}^\prime(u)\dd u+(-1)^r\widehat{\Phi}(0)\int_1^{+\infty}\frac{\theta(t)-t}{t^2}\dd t \\
+O\left(\frac{1}{\log^2{\left(q^r\right)}}\right)
\end{multline*}
and finally
\begin{equation*}
S_3=(-1)^r\frac{\log{\left(q^r\right)}}{4}\Phi(0)+(-1)^r\widehat{\Phi}(0)\left(1+\int_1^{+\infty}\frac{\theta(t)-t}{t^2}\dd t\right)+O\left(\frac{1}{\log^2{\left(q^r\right)}}\right).
\end{equation*}
We finally take care of the \underline{fourth term of \eqref{etape2}}. According to \eqref{eq_ams} and \eqref{eq_S}, we have %
\[%
 S\left(\alpha_f(p);n,r\right)=\delta_{2\mid r}+\sum_{\substack{1\leq j\leq r\\ j\equiv r\pmod{2}}}\left[\lambda_f(p^{jn})-\lambda_f(p^{jn-2})\right].
\]
One may remark that %
\[%
\sum_{\substack{p\in\mathcal{P} \\ p\nmid q}}\sum_{n\geq 3}\frac{\log{p}}{p^{n/2}}\widehat{\Phi}\left(\frac{n\log{p}}{\log{\left(q^r\right)}}\right)=\sum_{p\in\mathcal{P}}\sum_{n \geq 3}\frac{\log{p}}{p^{n/2}}\widehat{\Phi}\left(0\right)+O\left(\frac{1}{\log^3{(q^r)}}\right)
\]
since $\widehat{\Phi}(u)=\widehat{\Phi}(0)+O(u^2)$. Then, we easily get
\begin{equation*}
\sum_{p\in\mathcal{P}}\sum_{n\geq 3}\frac{\log{p}}{p^{n/2}}=\sum_{p\in\mathcal{P}}\frac{\log{p}}{p^{3/2}-p}.
\end{equation*}
\end{proof}

\section{Proof of Theorem \ref{thm_B}}

The aim of this part is to determine an asymptotic expansion of %
\begin{equation*}
\sum_{f\in\prim{\kappa}{q}}\omega_q(f)D_{1,q}[\Phi;r](f)=\Eh[q]\left(D_{1,q}[\Phi;r]\right). %
\end{equation*}
According to proposition \ref{explicit} and the proof of \cite[eq. (4.6) and (4.7)]{RiRo}, if %
\begin{equation}\label{eq_nupetit}%
\nu<\left(1-\frac{1}{2(\kappa-2\theta)}\right)\frac{2}{r^2}
\end{equation}
then %
\begin{multline}\label{explicitaverage}
\Eh[q]\left(D_{1,q}[\Phi;r]\right)=\left[\widehat{\Phi}(0)+\frac{(-1)^{r+1}}{2}\Phi(0)\right]+\frac{\widehat{\Phi}(0)}{\log{(q^r)}}\left[C_{\infty}+2(-1)^{r+1}C_{\text{PNT}}-2\delta_{2\mid r}C\right] \\ +\Eh[q]\left(P_q^3[\Phi;r](f)\right)+O\left(\frac{1}{\log^3{\left(q^r\right)}}\right).
\end{multline}
The first term in \eqref{explicitaverage} is the main term given in Theorem \ref{thm_B}. We now estimate the penultemate term of \eqref{explicitaverage} \emph{via} the trace formula given in \cite[Proposition 2.2]{RiRo}: %
\begin{equation}\label{eq_pu}%
\Eh[q]\left(P_q^3[\Phi;r]\right)=%
\PP{q,\mathrm{new}}3[\Phi;r]+\PP{q,\mathrm{old}}3[\Phi;r]
\end{equation}
where
\begin{equation*}
\PP{q,\mathrm{new}}3[\Phi;r]=-\frac{2}{\log{\left(q^r\right)}}\sum_{\substack{p\in\prem \\ p\nmid q}}%
\sum_{n\geq 3}\left[\sum_{\substack{1\leq j\leq r \\
j\equiv r\pmod{2}}}\left(\Delta_q(p^{jn},1)-\Delta_q(p^{jn-2},1)\right)\right]\frac{\log{p}}{p^{n/2}}\widehat{\Phi}\left(\frac{n\log{p}}{\log{\left(q^r\right)}}\right)
\end{equation*}
and
\begin{multline*}
\PP{q,\mathrm{old}}3[\Phi;r]=\frac{2}{q\log{\left(q^r\right)}}%
\sum_{\ell\mid q^\infty}\frac{1}{\ell}\sum_{\substack{p\in\prem \\ p\nmid q}}\sum_{n\geq 3}\left[\sum_{\substack{1\leq j\leq r \\
j\equiv r\pmod{2}}}\left(\Delta_1(p^{jn}\ell^2,1)-\Delta_1(p^{jn-2}\ell^2,1)\right)\right] \\
\frac{\log{p}}{p^{n/2}}\widehat{\Phi}\left(\frac{n\log{p}}{\log{\left(q^r\right)}}\right).
\end{multline*}
For $m\neq 1$ we have %
\begin{equation*}
\Delta_k(m,1)= 2\pi i^\kappa\sum_{\substack{c\geq 1 \\
k\mid c}}\frac{S(m,1;c)}{c}J_{\kappa-1}\left(\frac{4\pi\sqrt{m}}{c}\right)
\end{equation*}
where $S(m,1;c)$ is a Kloosterman sum. Let us estimate the new part which can be written as
\begin{equation*}
\PP{q,\mathrm{new}}3[\Phi;r]=%
-\frac{2%
(2\pi i^\kappa)
}{\log{\left(q^r\right)}}\sum_{\substack{1\leq j\leq r \\
j\equiv r\pmod{2}}}\sum_{n\geq 3}\left(\PP{q,\mathrm{new}}3[\Phi;r,jn]-\PP{q,\mathrm{new}}3[\Phi;r,jn-2]\right) %
\end{equation*}
where %
\begin{equation}\label{eq_inter}
\PP{q,\mathrm{new}}3[\Phi;r,k]=\sum_{\substack{p\in\prem\\ p\neq q}}\frac{\log p}{p^{n/2}}\widehat{\Phi}\left(\frac{\log{p}}{\log{\left(q^{r/n}\right)}}\right)\sum_{\substack{c\geq 1 \\ q\mid c}}\frac{S(p^k,1;c)}{c}J_{\kappa-1}\left(\frac{4\pi\sqrt{p^k}}{c}\right). %
\end{equation}
By \cite[lemma 3.10]{RiRo}, the $c$-sum in \eqref{eq_inter} is bounded by
\begin{equation*}
\frac{\tau(q)}{\sqrt{q}}\begin{cases}
\left(\frac{\sqrt{p^k}}{q}\right)^{1/2} & \text{if $p>q^{2/k}$,} \\
\left(\frac{\sqrt{p^k}}{q}\right)^{\kappa-1} & \text{otherwise.}
\end{cases}
\end{equation*}
We deduce %
\begin{align*}
 \sum_{n\geq 3}\PP{q,\mathrm{new}}3[\Phi;r,jn] &\ll\frac{\tau(q)}{q^{\kappa-1/2}}\sum_{n\geq 3}\sum_{p\leq q^{r\nu/n}}\frac{1}{p^{n/2}}p^{rn(\kappa-1)/2}\log p\\
&\ll \frac{\tau(q)}{q^{\kappa-1/2}}\sum_{3\leq n\leq\nu r\log{q}\left/\log{2}\right.} \frac{1}{n}q^{\nu r\left[\left((\kappa-1)r-1\right)n/2+1\right]/n} \\%
&\ll \frac{\tau(q)}{q^{\kappa-1/2}}q^{\nu r[(\kappa-1)r-1]/2}q^{\nu r/3}\log\log(3q)\\
&\ll \frac{1}{q^{1/2}}
\end{align*}
as soon as $\nu<2/r^2$ (and in particular if \eqref{eq_nupetit} is satisfied). We make the same computations for $jn-2$ and find then that $\PP{q,\mathrm{new}}3[\Phi;r,k]$ is an admissible error term. The old part is %
\[%
 \PP{q,\mathrm{old}}3[\Phi;r]=%
\frac{2%
(2\pi i^\kappa)
}{q\log{\left(q^r\right)}}\sum_{\substack{1\leq j\leq r \\
j\equiv r\pmod{2}}}\sum_{n\geq 3}\left(\PP{q,\mathrm{old}}3[\Phi;r,jn]-\PP{q,\mathrm{old}}3[\Phi;r,jn-2]\right) %
\]
where %
\[%
 \PP{q,\mathrm{old}}3[\Phi;r,k]=\sum_{\substack{p\in\prem \\ p\neq q}}\frac{\log p}{p^{n/2}}\widehat{\Phi}\left(\frac{\log p}{\log q^{r/n}}\right)\sum_{\ell\mid q^\infty}\frac{1}{\ell}\Delta_1(p^k\ell^2,1). %
\]
From \cite[eq (3.2) and (3.3)]{RiRo} we have %
\[%
 \sum_{\ell\mid q^\infty}\frac{1}{\ell}\Delta_1(p^k\ell^2,1)\leq 2(k+1) %
\]
so that %
\[%
 \sum_{n\geq 3}\PP{q,\mathrm{old}}3[\Phi;r,jn]\ll 1 %
\]
and similary for $\PP{q,\mathrm{old}}3[\Phi;r,jn-2]$. Finally $\Eh[q]\left(P_q^3[\Phi;r]\right)$ enters the $O(1/\log^3q^r)$ term. %
\appendix%

\section{Some comments on an aesthetic identity}
It is possible to prove on induction on $k_0\geq 1$ the following equality in $\Q[T]$: %
\begin{equation}\label{eq_pair2}
X_{2k_0}-X_{2k_0-2}=\sum_{j=0}^{k_0-1}\sum_{1\leq k_j<k_{j-1}<\cdots<k_1<k_0}(-1)^{j}\left[\prod_{i=0}^{j-1}\binom{2k_i}{k_i-k_{i+1}}\right]\left\{T^{2k_j}-\binom{2k_j} {k_j}\right\}.
\end{equation}
As a consequence, if $K\geq 1$ then %
\begin{equation}\label{eq_impair}
X_{2K+1}-X_{2K-1}=(-1)^{K}T\left(1+\sum_{1\leq k_0\leq K}(-1)^{k_0}X_{2k_0}-X_{2k_0-2}\right).
\end{equation}
Now, use \eqref{eq_lintch} with $r=1$ (so that $X_1=T$) to get from \eqref{eq_pair2} the equality %
\begin{multline*}
 X_{2k_0}-X_{2k_0-2}=\\\sum_{j=0}^{k_0-1}\sum_{1\leq k_j<k_{j-1}<\cdots<k_1<k_0}(-1)^{j}\left[\prod_{i=0}^{j-1}\binom{2k_i}{k_i-k_{i+1}}\right]\left[\sum_{\ell=0}^{2k_j}x(2k_j,1,\ell)X_\ell-\binom{2k_j}{k_j}X_0\right] %
\end{multline*}
and compare the coefficients of $X_0$ to obtain, thanks to \eqref{propriox} the equality %
\begin{equation*}
\sum_{j=0}^{k_0-1}\sum_{1\leq k_j<k_{j-1}<\cdots<k_1<k_0}(-1)^{j}\left[\prod_{i=0}^{j-1}\binom{2k_i}{k_i-k_{i+1}}\right]\binom{2k_j}{k_j}\frac{k_j}{1+k_j}=0. %
\end{equation*}
We could have expressed formulas \eqref{eq_pair2} and \eqref{eq_impair} in terms of Fourier coefficients of primitive forms to determine the lower order terms. However, this is definitely not the best way to proceed since it consists in decomposing the polynomial $X_K-X_{K-2}$ in the canonical basis of $\Q[T]$ and decomposing again each element of this canonical basis in the Chebyshev basis $(X_\ell)_{\ell\in\N}$. %

\section{S.J. Miller's identity and Chebychev polynomials}
S.J. Miller (\cite[Equation (3.12) Page 6]{Mil3}) recently proved that
\begin{equation}\label{eq_miller}
\alpha_f(p)^K+\beta_f(p)^K=\sum_{\substack{0\leq k \leq K \\
k\equiv K\pmod{2}}}c_{K,k}\lambda_f(p)^k
\end{equation}
where $c_{K,k}=0$ if $k\equiv K+1\pmod{2}$ and
\begin{align*}
c_{0,0} &= 0, \\
c_{2K,0} &= 2(-1)^K \quad(K\geq 1), \\
c_{2K,2L} &= \frac{2(-1)^{K+L}K(K+L-1)!}{(2L)!(K-L)!} \quad (1\leq L\leq K),\\
c_{2K+1,2L+1} &= \frac{(-1)^{K+L}(2K+1)(K+L)!}{(2L+1)!(K-L)!} \quad (0\leq L\leq K).
\end{align*}
We would like to give a quick proof of this identity, the crucial tool being Chebychev polynomials.
\begin{proof}[\proofname{} of equation \eqref{eq_miller}]
We know that
\begin{equation*}
\alpha_f(p)^K+\beta_f(p)^K=X_K(\lambda_f(p))-X_{K-2}(\lambda_f(p))
\end{equation*}
for $K\geq 2$ according to \eqref{eq_S4}. Thus, the proof consists in decomposing the polynomial $X_K-X_{K-2}$ in the canonical basis of $\mathbb{Q}[T]$. This can be done via \eqref{eq_lin}. It entails that
\begin{multline*}
\alpha_f(p)^K+\beta_f(p)^K=\sum_{\substack{0\leq u\leq K-2 \\
u\equiv K\pmod{2}}}(-1)^{(K-u)/2}\left[\binom{(K+u)/2}{u}+\binom{(K+u)/2-1}{u}\right]\lambda_f(p)^u \\
+\sum_{\substack{K-1\leq u\leq K \\
u\equiv K\pmod{2}}}(-1)^{(K-u)/2}\binom{(K+u)/2}{u}\lambda_f(p)^u,
\end{multline*}
which is an equivalent formulation of \eqref{eq_miller}.
\end{proof}
\begin{remark}
Equation \eqref{eq_miller} could be used to recover the lower order terms coming from $P_q^3[\Phi;r]$ but,once again, it is not the most clever way to proceed since it would imply decomposing the polynomials $X_K-X_{K-2}$ in the canonical basis of $\mathbb{Q}[T]$ at the beginning of the process and decomposing the polynomials $T^j$ in the basis $(X_r)_{r\geq 0}$ just before the end of the proof in order to be able to apply some trace formula for the Fourier coefficients of cusp forms.
\end{remark}

\providecommand{\bysame}{\leavevmode\hbox to3em{\hrulefill}\thinspace}
\providecommand{\MR}{\relax\ifhmode\unskip\space\fi MR }
\providecommand{\MRhref}[2]{%
  \href{http://www.ams.org/mathscinet-getitem?mr=#1}{#2}
}
\providecommand{\href}[2]{#2}

\end{document}